\documentclass{jnmp}

\begin{document}
\setcounter{page}{220}
\renewcommand{\evenhead}{P~Grozman and D~Leites}
\renewcommand{\oddhead}{The Shapovalov Determinant for the Poisson
Superalgebras}

\thispagestyle{empty}


\FirstPageHead{8}{2}{2001}
{\pageref{grozman-firstpage}--\pageref{grozman-lastpage}}{Article}

\copyrightnote{2001}{P~Grozman and D~Leites}

\Name{The Shapovalov Determinant 
 for the\\ Poisson Superalgebras}\label{grozman-firstpage}

\Author{Pavel GROZMAN and Dimitry LEITES}

\Address{Department of Mathematics, University of Stockholm\\
Roslagsv.  101, Kr\"aftriket hus 6, SE-106 91, Stockholm, Sweden\\
E-mail: mleites@matematik.su.se}

\Date{Received April 3, 2000; Revised January 24, 2001; Accepted
January 25, 2001}

\begin{abstract}
\noindent Among simple ${\mathbb Z}$-graded Lie superalgebras of
polynomial growth, there are several which have no Cartan matrix but,
nevertheless, have a quadratic even Casimir element $C_{2}$: these are
the Lie superalgebra $\mathfrak{k}^L(1|6)$ of vector fields on the
$(1|6)$-dimensional supercircle preserving the contact form, and the
series: the finite dimensional Lie superalgebra $\mathfrak{sh}(0|2k)$
of special Hamiltonian fields in $2k$ odd indeterminates, and the
Kac--Moody version of $\mathfrak{sh}(0|2k)$.  Using $C_{2}$ we compute
N.~Shapovalov determinant for $\mathfrak{k}^L(1|6)$ and
$\mathfrak{sh}(0|2k)$, and for the Poisson superalgebras
$\mathfrak{po}(0|2k)$ associated with $\mathfrak{sh}(0|2k)$. 
A.~Shapovalov described irreducible finite dimensional representations
of $\mathfrak{po}(0|n)$ and $\mathfrak{sh}(0|n)$; we generalize his
result for Verma modules: give criteria for irreducibility of the
Verma modules over $\mathfrak{po}(0|2k)$ and $\mathfrak{sh}(0|2k)$.
\end{abstract}

\section*{Introduction}

Every simple finite dimensional Lie algebra has a symmetrizable Cartan
matrix.  Moreover, if the simple ${\mathbb Z}$-graded Lie algebra $\mathfrak{g}$ of
polynomial growth (SZGLAPG, for short) has a Cartan matrix, i.e.,
$\mathfrak{g}=\mathfrak{g}(A)$, then $A$ is symmetrizable.  These Cartan matrices
correspond to Dynkin diagrams and extended Dynkin diagrams.  More
exactly, the algebras corresponding to extended diagrams are not
simple, they are certain ``relatives'', called Kac--Moody algebras, of
central extensions of simple ones; in applications Kac--Moody algebras
are even more interesting than simple ones, cf.~\cite{grozman:K3}.

For finite dimensional simple Lie algebras $\mathfrak{g}(A)$ N~Shapovalov
\cite{grozman:Sh} suggested a powerful method for description of irreducible
highest weight $\mathfrak{g}(A)$-modules (Verma modules and their quotients).
The method (a development of an idea due to Gelfand and Kirillov 
cf.~\cite{grozman:GK} with \cite{grozman:Sh}) was to consider what is now called the {\it
Shapovalov determinant}.

The Cartan matrices $A$ corresponding to finite dimensional Lie
algebras and to those of class SZGLAPG are very special.  Kac and
Kazhdan \cite{grozman:KK} extended Shapovalov's result to the Lie algebras
with {\it any} symmetrizable Cartan matrix $A$.  Namely, if $A$ is
symmetrizable, then $\mathfrak{g}(A)$ possesses a nondegenerate invariant
bilinear form $B$ and with the help of the associated with $B$
quadratic Casimir element $C_{2}$ they computed the Shapovalov
determinant for all such Lie algebras.

Using absence of zero divisors in the enveloping algebra of any Lie
algebra, they further obtained a description of irreducible modules
that occur in the Jordan--H\"older series of an arbitrary Verma module
over these algebras.

Kac later conjectured a formula for the Shapovalov determinant for the
Lie {\it super}algeb\-ras with symmetrizable Cartan matrix
\cite{grozman:K1}; the formula had an obvious mistake and a~correction
was offered in \cite{grozman:K2}, for the proof of the corrected
formula see \cite{grozman:M}.  As shown in \cite{grozman:KK}, the
formula for the determinant is a corollary of an explicit form of the
quadratic Casimir element.  Moreover, thanks to the existence of the
Casimir element, the Shapovalov determinant turns out to be equal to
the product of {\it linear} functions.  For Lie superalgebras with
symmetrizable Cartan matrix this element is even, so the argument
of~\cite{grozman:KK} are applicable literally.  Still, the direct
analogy soon stops: (a) to describe irreducible modules that occur in
the Jordan--H\"older series of an arbitrary Verma module over Lie
superalgebras we need new ideas due to the presence of zero divisors
(cf.~\cite{grozman:M}), (b) the form of the product of linear factors
depends on properties of odd roots involved and is more complicated,
see Th.~2.4~\cite{grozman:M}.  Still, one obtains a criterion for
irreducibility of Verma modules.

{\bf Our result}: calculation of the even quadratic Casimir element
for simple finite dimensional Lie superalgebras without symmetrizable
Cartan matrix.  Such are only the Lie superalgebras ${\mathfrak{sh}}(0|2n)$ of
special hamiltonian vector fields.  We also consider a ``relative'' of
${\mathfrak{sh}}(0|2n)$, the Poisson superalgebra $\mathfrak{po}(0|2n)$.  Actually we
derive the result for ${\mathfrak{sh}}(0|2n)$ from that for $\mathfrak{po}(0|2n)$.

Corollaries: a criterion for irreducibility of Verma modules over all
these algebras and calculation of the Shapovalov determinant: it is
the product of the linear terms-constituents of the above criterion.

Quantization sends $\mathfrak{sh}(0|2k)$ and $\mathfrak{po}(0|2k)$
into Lie superalgebras with Cartan mat\-rix~\cite{grozman:L1}; so our
result can be read as a ``dequantization'' of the Shapovalov
determinant for ${\mathfrak{psl}}\left(2^{k-1}|2^{k-1}\right)$ and
$\mathfrak{gl}\left(2^{k-1}|2^{k-1}\right)$.

Another corollary related with existence of the quadratic Casimir
element is an explicit solution to the classical Yang-Baxter equation
with values in the above Lie superalgebras; for finite dimensional
simple Lie superalgebras these solutions are given in~\cite{grozman:LSe}.

Since it is natural to describe Poisson Lie superalgebras as
subalgebras of the Lie superalgebra of contact vector fields, we
describe them and recall our earlier result on the Shapovalov
determinant for ${\mathfrak{k}}(1|6)$.

{\bf Related open problems}.  1) Extension of our result to loop
algebras is straightforward, elsewhere we will consider their
``Kac--Moody'' versions.

2) Even in the absence of the even quadratic Casimir element one can
define the Shapovalov determinant provided the algebra possesses an
involutive antiautomorphism which sends positive roots into negative
ones.  Among Lie superalgebras of type SZGLAPG only stringy or
``superconformal'' Lie superalgebras (for their complete list
see~\cite{grozman:GLS}) possess this property together with relatives
of the ``queer'' series and ${\mathfrak{sh}}(0|2n-1)$ as well as {\it
its} relative, $\mathfrak{po}(0|2n-1)$, together with loop algebras
and ``Kac--Moody'' versions thereof.  For the queer Lie superalgebras,
even finite dimensional ones, nobody had yet computed Shapovalov
determinant.  For some (but not all!)  of the {\it distinguished},
see~\cite{grozman:GLS}, stringy superalgebras the Shapovalov
determinant is computed.  In these cases it turned out to be the
product of indecomposable {\it quadratic} polynomials,
cf.~\cite{grozman:FF,grozman:KW,grozman:KR} and refs. therein.

3) After graded algebras it is natural to consider filtered ones,
prime examples being the Lie algebras of (a)~differential operators
with polynomial coefficients and of (b)~``complex size matrices''; for
a description of a large class of them see~\cite{grozman:GL3}; most of them
have nondegenerate invariant symmetric bilinear forms.  Though even
for the simplest of these algebras the Casimir element is not
calculated yet, Shoikhet~\cite{grozman:Sho} calculated the Shapovalov
determinant for some modules; conjecturally it is possible to
calculate Casimir element on a wider class of modules, cf.~\cite{grozman:LS}.

\section{A description of $\pbf{\mathfrak{k}^L (1|6)}$}

{\bf Supercircle.} A {\it supercircle} or (for a physicist) a
closed {\it superstring} of dimension $1|m$ is the real supermanifold
$S^{1|m}$ associated with the rank $m$ trivial vector bundle over the
circle.  Let $t=e^{i\varphi}$, where $\varphi$ is the angle parameter
on the circle, be the even indeterminate of the Fourier transforms;
let $\theta=(\theta_1, \dots, \theta_m)$, be the odd coordinates on
the supercircle formed by a basis of the fiber of the trivial bundle
over the circle.  Then $(t, \theta)$ are the coordinates on
$({\mathbb C}^*)^{1|m}$, the complexification of $S^{1|m}$.

Let $m=2k$ and the contact form be
\[
\alpha= dt -\sum_{1\leq i\leq k}(\xi_id\eta_i +\eta_id\xi_i) .
\]

For the classification of simple ``stringy" Lie superalgebras of
vector fields and their nontrivial central extensions see \cite{grozman:GLS}.
Among the ``main" series are: ${\mathfrak{vect}}^L(1|n)\!=\!{\mathfrak{der}}~{\mathbb C}[t^{-1}\!,
t,\theta]$, of all vector fields and $\mathfrak{k}^L(1|n)$ that preserves the
Pfaff equation $\alpha=0$.  The superscript $L$ indicates that we
consider vector fields with Laurent coefficients, not polynomial ones.

{\bf  The modules of tensor fields.} To advance further, we
have to recall the definition of the modules of tensor fields over
the general vectorial Lie superalgebra ${\mathfrak{vect}}(m|n)$ and its
subalgebras, see~\cite{grozman:BL}.  Let $\mathfrak{g}={\mathfrak{vect}}(m|n)$ realized by
vector fields on the $m|n$-dimensional linear supermanifold
${\mathcal C}^{m|n}$ with coordinates $x=(u, \xi)$ with the standard
grading ($\deg x_{i}=1$ for any $i=1$, \dots, $n+m$) and ${\mathfrak{g}}_{\geq 0}
=\mathop{\oplus}\limits_{i\geq
0}\mathfrak{g}_{i}$.  Clearly, $\mathfrak{g}_0\cong {\mathfrak{gl}}(m|n)$.  Let $V$ be the
${\mathfrak{gl}}(m|n)$-module with the {\it lowest} weight $\lambda=\mbox{lwt}(V)$.
Make $V$ into a ${\mathfrak{g}}_{\geq 0}$-module by setting ${\mathfrak{g}}_{+}\cdot V=0$
for ${\mathfrak{g}}_{+}=\mathop{\oplus}\limits_{i> 0}{\mathfrak{g}}_{i}$.  The
superspace $T(V)={\mbox{Hom}}_{U({\mathfrak{g}}_{\geq 0})}(U({\mathfrak{g}}), V)$ is isomorphic,
due to the Poincar\'e--Birkhoff--Witt theorem, to
${{\mathbb C}}[[x]]\otimes V$.  Its elements have a natural
interpretation as formal {\it tensor fields of type} $V$ (or
$\lambda$).  When $\lambda=(a, \dots , a)$ we will simply write
$T(\vec a)$ instead of $T(\lambda)$.  In what follows we consider
irreducible ${\mathfrak{g}}_0$-modules; for any other ${\mathbb Z}$-graded vectorial Lie
superlagebra construction of modules with lowest weight is identical.

{\bf Examples.} $T(\vec 0)$ is the superspace of functions;
$\mbox{Vol}(m|n)=T(1, \dots , 1; -1, \dots , -1)$ (the semicolon separates
the first $m$ coordinates of the weight with respect to the matrix
units $E_{ii}$ of ${\mathfrak{gl}}(m|n)$) is the superspace of {\it densities} or
{\it volume forms}.  We denote the generator of $\mbox{Vol}(m|n)$
corresponding to the ordered set of indeterminates $x$ by $\mbox{vol}(x)$.
The space of $\lambda$-densities is ${\mbox{Vol}}^{\lambda}(m|n)=T(\lambda,
\dots , \lambda; -\lambda, \dots , -\lambda)$.  In particular,
${\mbox{Vol}}^{\lambda}(m|0)=T(\vec \lambda)$ while
${\mbox{Vol}}^{\lambda}(0|n)=T(\overrightarrow{-\lambda})$.

{\bf Modules of tensor fields over stringy superalgebras.}
Denote by $T^L(V)= {\mathbb C}[t^{-1}, t]\otimes V$ the ${\mathfrak{vect}}(1|n)$-module
that differs from $T(V)$ by allowing the {\it Laurent} polynomials as
coefficients of its elements instead of just polynomials.  Clearly,
$T^L(V)$ is a ${\mathfrak{vect}}^L(1|n)$-module.  Define the {\it twisted with
weight $\mu$} version of $T^L(V)$ by
setting:
\[
T^L_\mu (V)={\mathbb C}[t^{-1}, t]t^\mu\otimes V.
\]

{\bf The ``simplest" modules}.  These are analogs of the {\it
standard} or {\it identity} representation of the matrix algebras.
The simplest modules over the Lie superalgebras of series ${\mathfrak{vect}}$
are, clearly, the modules ${\mbox{Vol}}^{\lambda}$.  These modules are
characterized by the fact that over~${\mathcal F}$, the algebra of functions,
they are of rank~1, i.e., have only one generator.  Over stringy
superalgebras, we can as well twist these modules and consider
${\mbox{Vol}}^\lambda_\mu$.  Observe that for $\mu\not\in{\mathbb Z}$ this module has
only one submodule, the image of the exterior differential $d$, see~\cite{grozman:BL};
for $\mu\in{\mathbb Z}$ this submodule coincides with the kernel
of the residue:
\[\arraycolsep=0em
\begin{array}{l}
\mbox{Res}: {\mbox{Vol}}^L \longrightarrow {\mathbb C},
\vspace{1mm}\\
 \displaystyle f\mbox{vol}(t, \xi) \mapsto
\ \mbox{the coeff.  of}\ \frac{\xi_{1} \ldots \xi_{n}}{t} \ \mbox{in
the power series expansion of}\ f.
\end{array}
\]
Over contact superalgebras $\mathfrak{k}(2n+1|m)$, it is more natural
to express the simplest modules not in terms of $\lambda$-densities
but in terms of powers of $\alpha$:
\[
{\mathcal F}_\lambda=\left\{ \begin{array}{ll}
{\mathcal F}\alpha^\lambda \quad  & \mbox{for} \ n=m=0,
\vspace{1mm}\\
{\mathcal F}\alpha^{\lambda/2} & \mbox{otherwise}.\end{array}\right.
\]
Observe that, as $\mathfrak{k}(2n+1|m)$-modules,
${\mbox{Vol}}^\lambda\cong{\mathcal F}_{\lambda(2n+2-m)}$
and ${\mathcal F}={\mathcal F}_{0}$.  In particular,
the Lie superalgebra of series $\mathfrak{k}$ does not distinguish between
$\frac{\partial}{\partial t}$ and $\alpha^{-1}$: their transformation
rules are identical.  Hence, $\mathfrak{k}(2n+1|m)\cong {\mathcal F}_{-1}$ if $n=m=0$ or
${\mathcal F}_{-2}$ otherwise.  (Physicists usually set $\deg \theta=\frac12$
and $\deg t=1$, whereas we prefer, as is customary among
mathematicians, the integer values of the highest weight with respect
to the Cartan subalgebra of $\mathfrak{k}(1|n)_{0}\cong {\mathfrak{o}}(2n)\oplus{\mathbb C} z$,
so we use doubled physicists weights.)

{\bf Convenient formulas.} A laconic way to describe the Lie
superalgebras of series $\mathfrak{k}$ is via {\it generating functions}.  For
$f\in{\mathbb C} [t, \theta]$, where $\theta=(\xi, \eta)$, set:
\[
K_f=(2-E)(f)\frac{\partial}{\partial t}-H_f + \frac{\partial f}{\partial t} E,
\qquad \mbox{where} \quad E=\sum\limits_i \theta_i
\frac{\partial}{\partial \theta_i}
\]
and where $H_f$ is the
hamiltonian field with Hamiltonian $f$ that preserves $d\alpha$:
\[
H_f=-(-1)^{p(f)}\sum\limits_{j\leq k}
\left(\frac{\partial f}{\partial \xi_j} \frac{\partial}{\partial \eta_j}+
\frac{\partial f}{\partial \eta_j} \frac{\partial}{\partial \xi_j}\right).
\]
Since
\begin{equation}
L_{K_f}(\alpha)=K_1(f)\alpha,
\end{equation}
it follows that $K_f\in \mathfrak{k} (2n+1|m)$.

To the (super)commutator $[K_f, K_g]$ there corresponds the {\it
contact bracket} of the
generating functions:
\[
[K_f, K_g]=K_{\{f, g\}_{k.b.}}.
\]
An explicit formula for the contact brackets is as follows.  Let us
first define the brackets on functions that do not depend on $t$.  The
{\it Poisson bracket} $\{\cdot , \cdot\}_{P.b.}$ is given by the
formula
\[
\{f, g\}_{P.b.}=-(-1)^{p(f)}\left [\sum\limits_{j\leq m}
\left(\frac{\partial f}{\partial \xi_j} \frac{\partial g}{\partial \eta_j}+
\frac{\partial f}{\partial \eta_j} \frac{\partial g}{\partial \xi_j}\right)\right].
\]
Now, the contact bracket is
\[
\{ f, g\}_{k.b.}=(2-E)(f)\frac{\partial g}{\partial t}-\frac{\partial f}{\partial t}(2-E)(g)-\{ f, g\}_{P.b.}.
\]
It is not difficult to prove that $\mathfrak{k} (1|2k)\cong\mbox{Span}\,(K_f: f\in
{\mathbb C}[t, \theta])$ as superspaces.

The {\it Poisson superalgebra} is $\mathfrak{po}(0|m)=\mbox{Span}\,(K_f: f\in
{\mathbb C}[\theta])$.  Its quotient modulo the center, ${\mathfrak{z}}={\mathbb C} K_1$, is
called the {\it Hamiltonian Lie superalgebra} ${\mathfrak{h}}(0|m)$; clearly,
${\mathfrak{h}}(0|m)\cong \mbox{Span}\,(H_f: f\in {\mathbb C}[p, q, \theta])$.
On $\mathfrak{po}(0|m)$ and ${\mathfrak{h}}(0|m)$, there are supertraces:
\[
K_f, \quad H_f\mapsto \int f\cdot \mbox{vol}(\theta).
\]
The traceless elements span ideals $\mathfrak{spo}(0|m)$, {\it special
Poisson superalgebra}, and ${\mathfrak{sh}}(0|m)$, {\it special Hamiltonian
superalgebra}.  For $m=2k$ quantization sends them into
$\mathfrak{gl}\left(2^{k-1}|2^{k-1}\right)$ and $\mathfrak{psl}\left(2^{k-1}|2^{k-1}\right)$, respectively,
and the integral becomes the usual supertrace~\cite{grozman:L1}.

{\bf Roots of $\pbf{\mathfrak{k}^L (1|6)}$.}
The Cartan subalgebra of ${\mathfrak{g}}=\mathfrak{k}^L (1|6)$ is the span of
\[\arraycolsep=0em
\begin{array}{llll}
H_1 = K_{\xi_1\eta_1}, & H_2 = K_{\xi_2\eta_2}, & H_3 = K_{\xi_3\eta_3}, & H_4 = K_{t},
\vspace{1mm}\\
H_5 = K_{\frac 1t\xi_2\xi_3\eta_3\eta_2}, \quad & H_6 = K_{\frac
1t\xi_1\xi_3\eta_3\eta_1},\quad  & H_7 = K_{\frac
1t\xi_1\xi_2\eta_2\eta_1}, \quad & H_8 = K_{\frac
1{t^2}\xi_1\xi_2\xi_3\eta_3\eta_2\eta_1}.
\end{array}
\]
The weight of the elements of ${\mathfrak{g}}$ is taken, however, only with
respect to the {\it diagonalizing} part ${\mathfrak{d}}$ of Cartan subalgebra,
namely, with respect to $H_4$ and, after semicolon, $H_1$, $H_2$,~$H_3$:
\[
\mbox{wht}(K_{t^a\xi_1^{\alpha_{1}}\xi_2^{\alpha_{2}}
\xi_3^{\alpha_{3}}\eta_3^{\beta_{3}}\eta_2^{\beta_{2}}
\eta_1^{\beta_{1}}})=(2(a-1)+|\alpha|+|\beta|; \alpha_{1}-\beta_{1},
\alpha_{2}-\beta_{2}, \alpha_{3}-\beta_{3}),
\]
where $|\alpha|=\sum \alpha_{i}$, $|\beta|=\sum \beta_{i}$.  The root
vectors are ordered lexicographically.

{\bf Roots of $\pbf{\mathfrak{po}(0|2k)}$ and $\pbf{{\mathfrak{sh}}(0|2k)}$.} The Cartan
subalgebra of $\mathfrak{po}(0|2k)$ is the span of $H_\emptyset = K_{1}$, $H_i
= K_{\xi_i\eta_i}$, $H_{i, j} = K_{\xi_i\xi_j\eta_i\eta_j}$, \dots,
$H_{1, \dots , n} = K_{\xi_1\dots \xi_k\eta_1\dots\eta_k}$.  In what
follows we denote: ${\mathcal I}=(1, \dots , n)$ and let ${\mathcal P}({\mathcal I})$ be the set of
all the subsets of ${\mathcal I}$.

Observe that the diagonalizing subalgebra ${\mathfrak{d}}$ of the Cartan
subalgebra of $\mathfrak{po}(0|2k)$ is the Cartan subalgebra of ${\mathfrak{o}}(2k)$
spanned by the $H_i$.  We will call it the ``small'' Cartan
subalgebra.  The roots are given with respect to ${\mathfrak{d}}$ completed with
(to enable a finer grading) the exterior derivation
$E=\sum\theta_i\frac{\partial}{\partial \theta_i}$, where $\theta=(\xi, \eta)$.  The
weight with respect to $E$, separated by a semicolon, is given first
(shifted by $-2$).

The root vectors and their weights are:
\[
\mbox{wht}(K_{\xi_1^{\alpha_{1}}\dots
\xi_k^{\alpha_{k}}\eta_1^{\beta_{1}}\dots \eta_k^{\beta_{k}}})=
(|\alpha|+|\beta|-2; \alpha_{1}-\beta_{1}, \dots,
\alpha_{k}-\beta_{k})
\]
and the root vectors are ordered lexicographically.
A shorthand for the vector \linebreak $K_{\xi_1^{\alpha_{1}}\dots
\xi_k^{\alpha_{k}}\eta_1^{\beta_{1}}\dots \eta_k^{\beta_{k}}}$ is
$K_{\xi^{\alpha}\eta^{\beta}}$.

The roots of ${\mathfrak{sh}}(0|2k)$ are similar to those of
$\mathfrak{po}(0|2k)$: just delete $K_{1}$ and $K_{\xi_1\eta_1\dots
\xi_k\eta_k}$ and replace all $K_{f}$ with $H_{f}$.

The root vectors corresponding to simple roots are
\[
X^+_{i; I}=K_{\xi_i\cdot \xi^I\eta^I} \quad \mbox{and} \quad
X^-_{i; I}=K_{\eta_i\cdot \xi^I\eta^I}
\quad \mbox{for any subset $I$ of ${\mathcal I}$ except ${\mathcal I}$ itself.}
\]
Though this does not matter in this paper, observe that these root
vectors $X^{\pm}_{i; I}$ are not the most natural generators of
$\mathfrak{po}(0|2k)$; at least, a seemingly more natural and more economic set
of generators is different, cf.~\cite{grozman:GLP}.

{\bf Kac--Moody superalgebras associated with $\pbf{\mathfrak{po}(0|2k)}$ and
$\pbf{{\mathfrak{sh}}(0|2k)}$.} Recall that for each finite dimensional Lie
superalgebra ${\mathfrak{g}}$, the {\it loop algebra} associated with it is
denoted by ${\mathfrak{g}}^{(1)}={\mathbb C}[t^{-1}, t]\otimes {\mathfrak{g}}$.  Let $B$ be the
nondegenerate supersymmetric bilinear form on ${\mathfrak{g}}$.  A~nontrivial
central extension $\widehat{{\mathfrak{g}}^{(1)}}={\mathfrak{g}}^{(1)}\oplus{\mathbb C} z$ of
${\mathfrak{g}}^{(1)}$ is given by means of the cocycle
\[
c(f, g) =\mbox{Res}\; B(f, dg)\quad \mbox{for any}\quad f, g\in {\mathfrak{g}}^{(1)}.
\]
On $\widehat{{\mathfrak{g}}^{(1)}}$, the form $B$ induces the following
nondegenerate invariant supersymmetric form $B^{(1)}$:
\[
B^{(1)}((f, a), (g, b))= B(f, g)(0) + ab\quad\mbox{for any $f, g\in
{\mathfrak{g}}^{(1)}$ and $a, b\in{\mathbb C}$}.
\]
The algebra $\widehat{{\mathfrak{g}}^{(1)}}\oplus{\mathbb C} t\frac{d}{dt}$ is called
(affine) {\it Kac--Moody} Lie superalgebra, cf.~\cite{grozman:K3}; for the
list of simple Kac--Moody Lie superalgebras see~\cite{grozman:FLS}.

The Cartan subalgebra of $\widehat{\mathfrak{po}(0|2k)^{(1)}}$ (resp.
$\widehat{{\mathfrak{sh}}(0|2k)^{(1)}}$) is the span of the Cartan subalgebra
of $\mathfrak{po}(0|2k)$ (resp.  ${\mathfrak{sh}}(0|2k)$) and the central element~$z$;
hence, the roots are given with respect to $t\frac{d}{dt}$ and
$E=\sum\theta_{i}\frac{\partial}{\partial \theta_{i}}$, and, after semicolon, the
``small'' Cartan subalgebra, spanned by the $H_{i}$.

The weights of the root vectors are
\[
\mbox{wht}(t^nK_{\xi^{\alpha}\eta^{\beta}})=(n;
|\alpha|+|\beta|-2; \alpha_{1}-\beta_{1}, \dots,
\alpha_{k}-\beta_{k})
\]
and the root vectors are ordered lexicographically.

\section{The bilinear forms and related Casimir elements}

In what follows the root elements are normed so that $B(e_\alpha,
e_\alpha^*) = 1$ for the series $\mathfrak{po}$ and ${\mathfrak{sh}}$ as well as for
$\mathfrak{k}^L (1|6)$.  The vector $e_\alpha^*$ is called the {\it right dual}
of $e_\alpha$, cf.~\cite{grozman:L2}.

In the realization of $\mathfrak{k}^L (1|6)$ by means of generating
functions the invariant nondegenerate even supersymmetric bilinear
form $B$ is given by the formula
\begin{equation}
B(K_f, K_g)=\mbox{Res}\; fg, \quad \mbox{where} \ \mbox{Res}\, (f) = \mbox{the coefficient
of} \quad \frac{\xi_1\dots\xi_3\eta_1\dots\eta_3}{t}.
\end{equation}
It is easy to verify directly that
\begin{equation}
H_1^* = H_5, \qquad H_2^* = H_6,\qquad H_3^* = H_7,\qquad H_4^* = H_8.
\end{equation}

\medskip

\noindent
{\bf Lemma 1.} {\it The following Casimir element corresponds to the
form $B$ and, therefore, belongs to (the completion of) the center of
the enveloping algebra of $\mathfrak{k}^L (1|6)$:
\begin{equation}
C_2 =
\sum_{\alpha>0} e_\alpha^* e_\alpha + \sum_{i=1}^4 H_iH_i^* +
4 H_5 + 2 H_6 - 4 H_8.
\end{equation}}

Observe, that if in formula (4) we replace the right dual elements
with the left dual ones, i.e., such that $B(e_\alpha^*, e_\alpha) =
1$, we obtain an element which does not belong to the center.

In the realization of $\mathfrak{po} (0|2k)$ by means of generating
functions, the invariant nondegenerate bilinear form $B$ is given by
the formula
\[
B(K_f, K_g)=\int fg \, \mbox{vol}\,(\xi, \eta).
\]
Clearly, this form induces an invariant nondegenerate form $B(H_f,
H_g)=B(K_f, K_g)$ on ${\mathfrak{sh}}(0|2k)$.  Obviously, up to a sign,
$H_i^*=H_{{\mathcal I}\setminus\{i\}}$, $H_{i, j}^*=H_{{\mathcal I}\setminus\{i, j\}}$,
etc.

\medskip

\noindent
{\bf Lemma 2.} {\it The following Casimir element corresponds to the
form $B$ and, therefore, belongs to the center of
the enveloping algebra of $\mathfrak{po} (0|2k)$:
\begin{equation}
C_2 = 2\mathop{\sum}\limits_{\alpha>0} e_\alpha^* e_\alpha +
\mathop{\sum}\limits_{J\in{\mathcal P}({\mathcal I})} H_JH_J^* +
(-2)^{k-1}H_{{\mathcal I}\setminus\{k\}}.
\end{equation}

The Casimir element for ${\mathfrak{sh}}(0|2k)$ is obtained from the above one
by deleting terms with $H_{\emptyset}$ and $ H_{{\mathcal I}}$.}

\medskip

Observe that the third summand in (5) is precisely the element of
the small Cartan subalgebra corresponding to the weight
``$2\rho$'', which is defined for any finite dimentional Lie
superalgebra as the halfsum of positive even roots minus the
halfsum of positive odd roots. It is remarkable that its form is
so simple.

In the realization of $\widehat{\mathfrak{po} (0|2k)}^{(1)}$ with
the help of generating functions the invariant nondegenerate bilinear
form $B$ is given by the formula
\[
B\left ((t^nK_f, a), (t^mK_g, b)\right )=\delta_{m+n, 0}\int
fg\,\mbox{vol}\, (\xi, \eta)+ab.
\]
Clearly, this form induces an invariant nondegenerate form on
$\widehat{{\mathfrak{sh}} (0|2k)}^{(1)}$.

The  Casimir element corresponding to the
invariant form $B$ belongs to the center of
the (completed) enveloping algebra of $\widehat{\mathfrak{po} (0|2k)}^{(1)}$.

{\bf On proofs.} The proof of Lemmas~1 and~2 is a direct
verification: it suffices to apply root vectors corresponding to
simple negative roots.  This is a routine done by hand, but tests are
much easier to perform with the help of Grozman's SuperLie package,
see~\cite{grozman:GL1}.

\section{The Shapovalov determinant.  Irreducible
Verma modules over $\pbf{\mathfrak{k}^L (1|6)}$,
and $\pbf{\mathfrak{po}(0|2k)}$ with its relatives}

Let $a= (a_1, \dots , a_m)$ be the highest weight of the Verma module
$M^a$ over ${\mathfrak{g}}$ and $b=(b_1, \dots , b_m)$ (here $m=8$ for $\mathfrak{k}^L
(1|6)$, and $m=2^{k-1}$ for $\mathfrak{po}(0|2k)$ and $m=2^{k-1}-2$ for
${\mathfrak{sh}}(0|2k)$) be one of the weights of $M^a$, i.e., $b=a-\sum
n_ir_i$, where the $r_i$ are positive roots and $n_i\in {\mathbb N}$ if the
root $r_i$ is even, $n_i=0$ or 1 if the root $r_i$ is odd, i.e., $\sum
n_ir_i$ is a quasiroot.  By applying $C_2$ to the highest weight
vector of the Verma module we get, as in~\cite{grozman:KK}, the following
theorems.

\medskip

\noindent
{\bf Theorem 1.} {\it Let ${\mathfrak{g}}=\mathfrak{k}^L (1|6)$.  The module $M^a$ is
irreducible if and only if for every quasiroot
\[\arraycolsep=0em
\begin{array}{l}
a_1a_5+a_2a_6+a_3a_7+a_4a_8\neq
a_1b_5+a_2b_6+a_3b_7+a_4b_8+b_1a_5+b_2a_6
\vspace{1mm}\\
\phantom{a_1a_5} {}+b_3a_7+b_4a_8 +b_1b_5+b_2b_6+b_3b_7+b_4b_8-4b_5-2b_6+4b_8,
\end{array}
\]
or, in other words, if and only if $2(a+\rho, \beta)\neq (\beta, \beta)$
for every quasiroot $\beta$ and $\rho= 2H_{5}+H_6-2H_8$.}

\medskip

\noindent
{\bf Theorem~2.} {\it Let ${\mathfrak{g}}=\mathfrak{po}(0|2k)$.  The module $M^a$ is
irreducible if and only if
\[
2(a+\rho, \beta)\neq (\beta, \beta)\quad
\mbox{for every quasiroot $\beta$ and $\rho= (-2)^{k-2}H_{{\mathcal I}\setminus\{k\}}$.}
\]
More explicitly, the above formula can be written as
\[
\sum\limits_{J\in{\mathcal P}({\mathcal I})}a_{J}a_{J^{*}}\neq
\sum\limits_{J\in{\mathcal P}({\mathcal I})}a_{J}b_{J^{*}}+
\sum\limits_{J\in{\mathcal P}({\mathcal I})}b_{J}a_{J^{*}}+
\sum\limits_{J\in{\mathcal P}({\mathcal I})}b_{J}b_{J^{*}}-(-2)^{k-1}
b_{{\mathcal I}\setminus\{k\}}.
\]

The description of the irreducible Verma modules $M^a$ over
${\mathfrak{sh}}(0|2k)$ is similar: in the above formula delete the terms with
subscripts $\emptyset$ and ${\mathcal I}$.}

\subsection*{Acknowledgements}

The financial support of PG by the Swedish Institute and
same of DL by NFR and by the Hebrew University  of Jerusalem during
Forschheimer Professorship is most gratefully acknowledged.

\label{grozman-lastpage}

\end{document}